\documentclass[11pt]{article}
\usepackage[fleqn]{amsmath}
\usepackage{amsfonts}
\usepackage{amssymb}
\usepackage{hyperref}
\usepackage{indentfirst}
\usepackage{mathrsfs}
\usepackage{tikz}
\usepackage{mathtools}
\usepackage{graphicx}
\usepackage{xparse}

\usepackage[top=1in,bottom=1in,left=1.25in,right=1.25in]{geometry}
\date{}
\allowdisplaybreaks

\begin{document}

 \renewcommand{\baselinestretch}{1.2}
 \renewcommand{\arraystretch}{1.0}

 \title{\bf Yetter-Drinfeld modules for group-cograded Hopf quasigroups}
 \author
 { \textbf{Huili Liu} \footnote{College of Science, Nanjing Agricultural University, Nanjing 210095, Jiangsu, China.
  E-mail: 2020111006@njau.edu.cn}, \,
  \textbf{Tao Yang} \footnote{Corresponding author. College of Science, Nanjing Agricultural University, Nanjing 210095, Jiangsu, China.
             E-mail: tao.yang@njau.edu.cn}, \,
  \textbf{Lingli Zhu} \footnote{College of Science, Nanjing Agricultural University, Nanjing 210095, Jiangsu, China. E-mail: zhounan0805@163.com}
 }
 \maketitle

 \begin{center}
 \begin{minipage}{12.cm}
 {Abstract: Let $H$ be a crossed group-cograded Hopf quasigroup. We first introduce the notion of $p$-Yetter-Drinfeld quasimodule over $H$. If the antipode of $H$ is bijective,  we show that the category $\mathscr Y\mathscr D\mathscr Q(H)$ of Yetter-Drinfeld quasimodules over $H$ is a crossed category, and the subcategory $\mathscr Y\mathscr D(H)$ of Yetter-Drinfeld modules is a braided crossed category.
 }

 \\

 \end{minipage}
 \end{center}
 \normalsize

 \section{Introduction}
 \def\theequation{\thesection.\arabic{equation}}
 \setcounter{equation}{0}
 Hopf quasigroups were introduced by Klim and Majid in \cite{KM10} in order to understand the structure and relevant properties of the algebraic 7-sphere. They are non-associative generalizations of Hopf algebras, however there are certain conditions about antipode that can compensate for their lack of associativity. Hopf quasigroups are no longer  associative algebras, so its compatibility conditions are quite different from those of Hopf algebras.

 Turaev introduced the notion of braided crossed categories in \cite{T} which is based on a group $G$ and showed that such a category gives rise to a 3-dimensional homotogy quantum field theory with target space $K(G, 1)$. In fact, braided crossed categories are braided monoidal categories in Freyd-Yetter categories of crossed G-sets (see \cite{FY89}), play a key role in the construction of these homotopy invariants.

 Zunino introduced a kind of Yetter-Drinfeld modules over crossed group coalgebra in \cite{ZU04}, and constructed a braided crossed category by this kind of Yetter-Drinfeld modules. This idea were generalized to multiplier Hopf T-coalgebras by Yang in \cite{YW11}.
 It is natural to propose a question: Does this method also hold for some other algebraic structures?

 Motivated by this question, the main purpose of this paper is to construct a braided crossed category by $p$-Yetter-Drinfeld modules over crossed group-cograded Hopf quasigroups.

 This paper is organized as follows: In section 2, we  recall some notions, such as braided crossed categories, Turave's left index notation, and Hopf quasigroups.
 All of these are the most important building blocks for completing this article.

 In section 3, we introduce crossed group-cograded Hopf quasigroups and then give some examples of this algebraic structure.
 Moreover, we give a method to construct crossed group-cograded Hopf quasigroups, which is rely on a fixed crossed group-cograded Hopf quasigroup.
 At the end of this section, we show that a group-cograded Hopf quasigroup with the group G is indeed a Hopf quasigroup in Turaev category.

 In section 4, we first give the definition of $p$-Yetter-Drinfeld quasimodules over a crossed group-cograded Hopf quasigroups H.
 And then we show the category $\mathscr Y\mathscr D\mathscr Q(H)$ of Yetter-Drinfeld quasimodules over $H$ is a crossed category, and the subcategory $\mathscr Y\mathscr D(H)$ of Yetter-Drinfeld modules is a braided crossed category.

 \section{Preliminaries}

 \subsection{Crossed categories and Turaev category}

 Recall from \cite{JS93, ML71, T} that
 let $G$ be a group and $Aut(\mathscr C)$ be the group of invertible strict tensor functors from $\mathscr C$ to itself. A category $\mathscr C$ over $G$ is called a crossed category if it satisfies the following:
 \begin{itemize}
   \item[(1)] $\mathscr C$ is a monoidal category;
   \item[(2)] $\mathscr C$ is disjoint union of a family of subcategories $(\mathscr C_{\alpha})_{\alpha\in G}$, and for any $U\in \mathscr C_{\alpha}, V\in \mathscr C_{\beta}, U\otimes V\in \mathscr C_{\alpha \beta}$. The subcategory $\mathscr C_{\alpha}$ is called the $\alpha th$ component of $\mathscr C$;
   \item[(3)] Consider a group homomorphism $\phi : G\rightarrow Aut(\mathscr C), \beta \rightarrow \phi_{\beta}$, and assume that $\phi_{\beta}(\mathscr C_{\alpha})=\mathscr C_{\beta \alpha\beta^{-1}}$, for all $\alpha,\beta\in G$. The functors $\phi_{\beta}$ are called conjugation isomorphisms.
 \end{itemize}

 We will use Turave's left index notation from \cite{T, V02} for functors $\phi_{\beta}$: Given $\beta \in G$ and an object $V\in \mathscr C$, the functor $\phi_{\beta}$ will be denoted by $^{\beta}(\cdot)$ or $^{V}(\cdot)$ and $^{\beta^{-1}}(\cdot)$ will be denoted by ${}^{\bar{V}}(\cdot)$. Since ${}^{V}(\cdot)$ is a functor, for any object $U\in \mathscr C$ and any composition of morphism $g\circ f$ in $\mathscr C$, we obtain ${}^{V}id_{U}=id_{{}^{V}U}$ and $^{V}(g\circ f)={}^{V}g\circ {}^{V}f$. Since the conjugation $\varphi: G\rightarrow Aut(\mathscr C)$ is a group homomorphism, for any $V, W\in \mathscr C$, we have ${}^{V\otimes W}(\cdot)={}^{V}({}^{W}(\cdot)$ and $^{1}(\cdot)={}^{V}(^{\bar{V}}(\cdot))={}^{\bar{V}}({}^{V}(\cdot))=id_{\mathscr C}$. Since for any $V\in \mathscr C$, the functor $^{V}(\cdot)$ is strict, we have ${}^{V}(g\otimes f)={}^{V}g\otimes {}^{V}f$ for any morphism $f$ and $g$ in $\mathscr C$, and ${}^{V}(1)=1$.
 \\

 A braiding of a crossed category $\mathscr C$ is a family of isomorphisms $(C= C_{U,V})_{U,V\in \mathscr C}$, where
 $C_{U,V}: U\otimes V\rightarrow {}^{U}V\otimes U$satisfying the following conditions:

 \begin{itemize}
   \item[(1)] For any arrow $f\in \mathscr C_{p}(U,U')$ and $g\in \mathscr C(V,V')$,
   \begin{eqnarray}
   (({}^{p}g)\otimes f)\, C_{U,V}=C_{U',V'}\, (f\otimes g);
   \end{eqnarray}
   \item[(2)] For all $U, V, W\in \mathscr C$, we have
   \begin{eqnarray}
   C_{U\otimes V, W} &=& a_{{}^{U\otimes V}W, U, V}\, (C_{U,{}^{V}W}\otimes id_{V})\, a^{-1}_{U, {}^{V}W, V}\, (\iota_{U}\otimes C_{V,W})\, a_{U,V,W},\\
   C_{U, V\otimes W} &=& a^{-1}_{{}^{U}V, {}^{U}W, U}\, (\iota_{{}^{U}V}\otimes C_{U, W})\, a_{{}^{U}V, U, W}\, (C_{U,V}\otimes \iota_{W})\, a^{-1}_{U,V,W},
   \end{eqnarray}
   where a is the natural isomorphisms in the tensor category $\mathscr C$.
   \item[(3)] For all $U, V\in \mathscr C$ and $q\in G$,
   \begin{eqnarray}
   \phi_{q}(C_{U,V})=C_{\phi_{q}(U),\phi_{q}(V)}.
   \end{eqnarray}
 \end{itemize}

 A crossed category endowed with a braiding is called a braided crossed category.
 \\

 Turaev category as an special symmetric monoidal category is introduced by Caenepeel from \cite{CD06}. We recall the notion of Turaev category $\mathscr T_{R}$:
 Let $R$ be a commutative ring. A Turaev $R$-module is a couple $\underline{M}=(X, M)$, where $X$ is a set, and $M=(M_{x})_{x\in X}$ is a family of $R$-modules indexed by $X$. A morphisms between two $T$-modules $(X,M)$ and $(Y,N)$ is a couple $\underline{\phi}=(f,\phi)$, where $f: Y\rightarrow X$ is a function, and $\phi=(\phi_{y}: M_{f_{y}}\rightarrow N_{y})_{y\in Y}$ is a family of linear maps indexed by $Y$. The composition of $\underline{\phi}: \underline{M}\rightarrow \underline{N}$ and $\underline{\psi}: \underline{N}\rightarrow \underline{P}=(Z,P)$ is defined as follows:
 \begin{eqnarray*}
 &&\underline{\psi}\, \underline{\phi}=(fg, (\psi_{z}\phi_{g(z)})_{z\in Z}).
 \end{eqnarray*}
 The category of Turaev $R$-modules is called Turaev category and denoted by $\mathscr{T}_{R}$.

\subsection{Hopf quasigroups}

 Throughout this article, all spaces we considered are over a fixed field $k$.

 Recall from \cite{KM10} that a Hopf quasigroup $H$ is a unital (not necessarily associative) algebra $(H, \mu_{H}, \eta_{H})$ and a counital and coassociative coalgebra $(H, \delta_{H}, \epsilon_{H})$ with the morphisms $\delta_{H}$ and $\epsilon_{H}$ are algebra morphisms. And there exists a linear map $S: H\rightarrow H$ such that
 \begin{eqnarray}
 \mu(id\otimes \mu)(S\otimes id\otimes id)(\Delta \otimes id)=\epsilon \otimes id=\mu (id\otimes \mu)(id\otimes S\otimes id)(\Delta \otimes id ), \label{HQ1}\\
 \mu(\mu \otimes id)(id\otimes id\otimes S)(id\otimes \Delta)=id\otimes \epsilon=\mu (\mu \otimes id)(id\otimes S\otimes id)(id \otimes \Delta). \label{HQ2}
 \end{eqnarray}

 In this paper we use Sweelder notation for coproduct: $\Delta(h)=h_{(1)}\otimes h_{(2)}$, for $h\in H$. Use this notation we can rewrite the conditions (\ref{HQ1}) and (\ref{HQ2}) of Hopf quasigroup as
 \begin{eqnarray}
 S(h_{(1)})(h_{(2)}g)=\epsilon(h)g=h_{(1)}(S(h_{(2)})g),\label{HQ3}\\
 (gh_{(1)})S(h_{(2)})=g\epsilon(h)=(gS(h_{(1)}))h_{(2)},\label{HQ4}
 \end{eqnarray}
 for all $h,g\in H$.

 If the antipode $S$ of $H$ is bijective, then for all $h,g\in H$, we have
 \begin{eqnarray}
 &&S^{-1}(h_{2})(h_{1}g)=\epsilon(h)g=h_{2}(S^{-1}(h_{1})g), \label{HQ5}\\
 &&(gS^{-1}(h_{2}))h_{1}=g\epsilon(h)=g(h_{2}S^{-1}(h_{1})). \label{HQ6}
 \end{eqnarray}

 A morphism between Hopf quasigroups $H$ and $B$ is a map $f: H\rightarrow B$ which is both algebra and coalgebra morphism. And a Hopf quasigroup is associative if and only if it is a Hopf algebra. More details see \cite{BJ12, KM10}.

\section{Crossed group-cograded Hopf quasigroups}
 \def\theequation{\thesection.\arabic{equation}}
 \setcounter{equation}{0}

 In this section, we first introduced the notion of crossed group-cograded Hopf quasigroup, generalizing crossed Hopf group-coalgebra introduced in \cite{T}.
 Then we prove that a group-cograded Hopf quasigroup is indeed a Hopf quasigroup in Turaev category, and provide a method to construct crossed group-cograded Hopf quasigroups.
 \\

 \textbf{Definition \thesection.1}
 Let $G$ be a group. $(H=\bigoplus_{p\in G}H_{p}, \Delta, \epsilon)$ is called a group-cograded Hopf quasigroup over $k$, where each $H_{p}$ is an unital $k$-algebra with multiplication $\mu_{p}$ and unit $\eta_{p}$, comultiplication $\Delta$ is a family of homomorphisms $(\Delta_{p,q}: H_{pq}\rightarrow H_{p}\otimes H_{q})_{p,q\in G}$ and counit $\epsilon$ is a homomorphism defined by $\epsilon: H_{e}\rightarrow k$ such that the following conditions:
 \begin{itemize}
  \item[(1)] $H_{p}H_{q}=0$ whenever $p,q \in G$ and $p\neq q$, and $\eta_{p}(1_{k}) = 1_{p}$;

  \item[(2)] $\Delta$ is coassociative in the sense that for any $p, q, r\in G$,
   \begin{eqnarray}
   (\Delta_{p,q}\otimes id_{H_{r}})\Delta_{pq,r}=(id_{H_{p}}\otimes \Delta_{q,r})\Delta_{p,qr}, \label{GHde}
   \end{eqnarray}
  and for all $p, q\in G$ the $\Delta_{p,q}$ is an algebra homomorphism and $\Delta_{p,q}(H_{pq})\subseteq H_{p}\otimes H_{q}$.

  \item[(3)] $\epsilon$ is counitary in the sense that for any $p\in G$,
   \begin{eqnarray}
   (id_{H_{p}}\otimes \epsilon)\Delta_{p,e}=(\epsilon\otimes id_{H_{p}})\Delta_{e,p}=id_{H_{p}}, \label{GHep}
   \end{eqnarray}
  and $\epsilon$ is an algebra homomorphism and $\epsilon(1_{e})=1_{k}$;

  \item[(4)] endowed $H$ with algebra anti-homomorphisms $S=(S_{p}: H_{p}\rightarrow H_{p^{-1}})_{p\in G}$, then for any $p\in G$,
  \begin{eqnarray}
  \epsilon\otimes id_{H_{p}}&=&\mu_{p}(id_{H_{p}}\otimes \mu_{p})(S_{p^{-1}}\otimes id_{H_{p}}\otimes id_{H_{p}})(\Delta_{p^{-1},p}\otimes id_{H_{p}})\nonumber \\
  &=&\mu_{p}(id_{H_{p}}\otimes \mu_{p})(id_{H_{p}}\otimes S_{p^{-1}}\otimes id_{H_{p}})(\Delta_{p,p^{-1}}\otimes id_{H_{p}}),\label{S1}\\
  id_{H_{p}}\otimes \epsilon&=&\mu_{p}(\mu_{p}\otimes id_{H_{p}})(id_{H_{p}}\otimes id_{H_{p}}\otimes S_{p^{-1}})(id_{H_{p}}\otimes \Delta_{p,p^{-1}})\nonumber \\
  &=&\mu_{p}(\mu_{p}\otimes id_{H_{p}})(id_{H_{p}}\otimes S_{p^{-1}}\otimes id_{H_{p}})(id_{H_{p}}\otimes \Delta_{p^{-1},p}).\label{S2}
  \end{eqnarray}
 \end{itemize}


 We extend the Sweedler notation for a comultiplication in the following way:
 For any $p,q\in G, h_{pq}\in H_{pq}$,
 \begin{eqnarray*}
 \Delta_{p,q}(h_{pq})=h_{(1,p)}\otimes h_{(2,q)}.
 \end{eqnarray*}
 Then, we can rewrite the conditions (\ref{S1}) and (\ref{S2})as: for all $p\in G$ and $h_{e}\in H_{e}, g\in H_{p}$,
 \begin{eqnarray}
 S_{p^{-1}}(h_{(1,p^{-1})})(h_{(2,p)}g)=\epsilon(h_{e})g=h_{(1,p)}(S_{p^{-1}}(h_{(2,p^{-1})})g),\label{S3}\\
 (gh_{(1,p)})S_{p^{-1}}(h_{(1,p^{-1})})=g\epsilon(h_{e})=(gS_{p^{-1}}(h_{(1,p^{-1})}))h_{(2,p)}.\label{S4}
 \end{eqnarray}
 \smallskip

 As in the Hopf group-coalgebra (or group-cograded Hopf algebra) case,  we show group-cograded Hopf quasigroups are Hopf quasigroups in a special category as follows.
 \\

 \textbf{Proposition \thesection.2} If $H=\bigoplus_{p\in G} H_{p}$ is a group-cograded Hopf quasigroup, then $(G,H)$ is a Hopf quasigroup in Turaev category $\mathscr{T}_{k}$.
 \smallskip

 \textbf{\emph{Proof}} As $H$ is a group-cograded Hopf quasigroup and $G$ is a group, then we can give $\underline{H}=(G, H)$ a unital algebra structure $(\underline{H}, \underline{\mu}, \underline{\eta})$ by
 \begin{eqnarray*}
 \begin{aligned}
 \underline{k}&\stackrel{\underline{\eta}}{\longrightarrow}\underline{H}\\
 (*)&\stackrel{e}{\longleftarrow} G\\
 k&\stackrel{\eta_{p}} {\longrightarrow}H_{p},
 \end{aligned}
 \qquad \rm{and} \qquad
 \begin{aligned}
 \underline{H}\otimes \underline{H}&\stackrel{\underline{\mu}}{\longrightarrow}\underline{H}\\
 G\times G&\stackrel{\delta}{\longleftarrow} G\\
 H_{p}\otimes H_{p}&\stackrel{\mu_{p}} {\longrightarrow}H_{p},
 \end{aligned}
 \end{eqnarray*}
 such that
 \begin{eqnarray*}
 \begin{aligned}
 \underline{H}&\stackrel{\underline{id}\otimes \underline{\eta}}{\makebox[1cm]{\rightarrowfill}}&\underline{H}\otimes \underline{H}&\stackrel{\underline{\mu}}{\makebox[1cm]{\rightarrowfill}}&\underline{H}\\
 G&\stackrel{(1,e)}{\makebox[1cm]{\leftarrowfill}}&G\times G&\stackrel{\delta}{\makebox[1cm]{\leftarrowfill}}&G\\
 H_{p}&\stackrel{id\otimes \eta_{p}}{\makebox[1cm]{\rightarrowfill}}&H_{p}\otimes H_{p}&\stackrel{\eta_{p}}{\makebox[1cm]{\rightarrowfill}}&H_{p},
 \end{aligned}
 \qquad \rm{and} \qquad
 \begin{aligned}
 \underline{H}&\stackrel{\underline{\eta}\otimes\underline{id}}{\makebox[1cm]{\rightarrowfill}}&\underline{H}\otimes \underline{H}&\stackrel{\underline{\mu}}{\makebox[1cm]{\rightarrowfill}}&\underline{H}\\
 G&\stackrel{(e,1)}{\makebox[1cm]{\leftarrowfill}}&G\times G &\stackrel{\delta}{\makebox[1cm]{\leftarrowfill}}&G\\
 H_{p}&\stackrel{\eta_{p}\otimes id}{\makebox[1cm]{\rightarrowfill}}&H_{p}\otimes H_{p}&\stackrel{\eta_{p}}{\makebox[1cm]{\rightarrowfill}}&H_{p}.
 \end{aligned}
 \end{eqnarray*}

 Also, we can give $(G, H)$ a coalgebra structure $(\underline{H},\underline{\Delta},\underline{\epsilon})$ by
 \begin{eqnarray*}
 \begin{aligned}
 \underline{H}&\stackrel{\underline{\epsilon}}{\longrightarrow}\underline{k}\\
 G&\stackrel{i}{\longleftarrow}(*)\\
 H_{1}=H_{i(e)}&\stackrel{\epsilon}{\longrightarrow}k,
 \end{aligned}
 \qquad \rm{and} \qquad
 \begin{aligned}
 \underline{H}&\stackrel{\underline{\Delta}}{\longrightarrow}\underline{H}\otimes \underline{H}\\
 G&\stackrel{m}{\longleftarrow}G\times G\\
 H_{gh}=H_{m(gh)}&\stackrel{\Delta_{g,h}}{\longrightarrow} H_{g}\otimes H_{h},
 \end{aligned}
 \end{eqnarray*}
 such that $(\underline{\Delta},\underline{\epsilon})$ are algebra maps.

 Let $s: G\rightarrow G, s(g)=g^{-1}$,  then we can consider a map $\underline{S}=(s,S)$ in Turaev category as the antipode of $\underline{H}$, where $S$ is the antipode of group-cograded Hopf quasigroup $H$.
 Next, we will only check that $\underline{S}$ satisfy the condition (\ref{HQ3}), the condition (\ref{HQ4}) is similar. Indeed,
 \begin{eqnarray*}
 \begin{aligned}
 \underline{H}\otimes \underline{H}&\stackrel{\underline{\Delta}\otimes \underline{id}}{\makebox[1cm]{\rightarrowfill}}&\underline{H}\otimes \underline{H}\otimes \underline{H}&\stackrel {\underline{S}\otimes \underline{id}\otimes \underline{id}}{\makebox[1.5cm]{\rightarrowfill}}&\underline{H}\otimes \underline{H}\otimes \underline{H}&\stackrel{\underline{id}\otimes \underline{\mu}}{\makebox[1cm]{\rightarrowfill}}&\underline{H}\otimes \underline{H}&\stackrel{\underline{\mu}}{\makebox[0.6cm]{\rightarrowfill}}&\underline{H}\\
 G\times G&\stackrel{(m,1)}{\makebox[1cm]{\leftarrowfill}}&G\times G\times G&\stackrel{(s,1,1)}{\makebox[1.5cm]{\leftarrowfill}}&G\times G\times G &\stackrel{(1,\delta)}{\makebox[1cm]{\leftarrowfill}}&G\times G&\stackrel{\delta}{\makebox[0.6cm]{\leftarrowfill}}&G\\
 H_{e}\otimes H_{p}&\stackrel{\Delta\otimes id}{\makebox[1cm]{\rightarrowfill}}&H_{p^{-1}}\otimes H_{p}\otimes H_{p}&\stackrel{S_{p^{-1}}\otimes id \otimes id}{\makebox[1.6cm]{\rightarrowfill}}&H_{p}\otimes H_{p}\otimes H_{p}&\stackrel{id\otimes \mu_{p}}{\makebox[1cm]{\rightarrowfill}}&H_{p}\otimes H_{p}&\stackrel{\mu_{p}}{\makebox[0.6cm]{\rightarrowfill}}&H_{p},
 \end{aligned}
 \end{eqnarray*}
 and
 \begin{eqnarray*}
 \underline{H}\otimes \underline{H}&\stackrel{\underline{\epsilon}\otimes \underline{id}}{\makebox[1.5cm]{\rightarrowfill}}&\underline{k}\otimes \underline{H}\\
 G\times G&\stackrel{i\otimes 1}{\makebox[1.5cm]{\leftarrowfill}}& (*)\times G\\
 H_{i(e)}\otimes H &\stackrel{\epsilon\otimes id}{\makebox[1.5cm]{\rightarrowfill}}& k\otimes H.
 \end{eqnarray*}

 Since $H$ is a group-cograded Hopf quasigroup, we have $(\underline{\Delta}\otimes \underline{id})(\underline{S} \otimes\underline{id}\otimes\underline{id})(\underline{id}\otimes\underline{\mu})\underline{\mu}=\underline{\epsilon}\otimes \underline{id}$. Thus the left hand of equation (\ref{HQ3}) holds, and the right hand is similar. $\hfill \square$
 \\

 \textbf{Definition \thesection.3} A group-cograded Hopf quasigroup $(H=\bigoplus_{p\in G} H_{p},\Delta, \epsilon, S)$ is said to be a crossed group-cograded Hopf quasigroup provided it is endowed with a crossing $\pi :G\rightarrow Aut(H)$ such that
 \begin{itemize}
  \item[(1)] each $\pi_{p}$ satisfies $\pi_{p}(H_{q})=H_{pqp^{-1}}$, and preserves the counit, the antipode, and the comultiplication, i.e., for all $p, q, r\in G$,
   \begin{eqnarray}
   &&\epsilon\pi_{p}|_{H_{e}}=\epsilon, \label{pie}\\
   &&\pi_{p}S_{q}=S_{pqp^{-1}}\pi_{p},\label{pis}\\
   &&(\pi_{p}\otimes\pi_{p})\Delta_{q,r}=\Delta_{pqp^{-1},prp^{-1}}\pi_{p}, \label{pid}
   \end{eqnarray}
  \item[(2)] $\pi$ is multiplicative in the sense that for all $p,q\in G$,
  $\pi_{pq}=\pi_{p}\pi_{q}$.
 \end{itemize}

  If all of its subalgebras $(H_{p})_{p\in G}$ are associative, then $H$ is a crossed Hopf group-coalgebra introduced in \cite{T}.
  In the following we give two examples of crossed group-cograded Hopf quasigroups, and both examples are derived from an action of $G$ on a Hopf quasigroup over $k$ by Hopf quasigroup endomorphisms.
 \\

 \textbf{Example \thesection.4}  Let $(H, \Delta, \epsilon, S)$ be a Hopf quasigroup. Set $H^{G}=\big(H_{p})_{p\in G}$ and $G$ is the homomorphism group of $H$, where for each $p \in G$, the algebra
 $H_{p}$ is a copy of $H$. Fix an identification isomorphism of algebras $i_{p}: H\rightarrow H_{p}$. For $p, q \in G$, we define a comultiplication
 $\Delta_{p, q}: H_{pq}\rightarrow H_{p}\otimes H_{q}$ by
 \begin{eqnarray*}
 &&\Delta_{p, q}(i_{pq}(h))=\sum_{(h)}i_{p}(h_{(1)})\otimes i_{q}(h_{(2)}),
 \end{eqnarray*}
 where $h\in H$. The counit $\epsilon: H_{e}\rightarrow k$ is defined by $\epsilon(i_{e}(h))=\epsilon(h)\in k$ for $h\in H$. For $p \in G$, the antipode $S_{p}: H_{p}\rightarrow H_{p^{-1}}$ is given by
 \begin{eqnarray*}
 && S_{p}(i_{p}(h))=i_{p^{-1}}(S(h)),
 \end{eqnarray*}
 where $h\in H$. For $p, q \in G$, the homomorphism $\pi_{p}: H_{q}\rightarrow H_{pqp^{-1}}$ is defined by $\pi_{p}(i_{q}(h))=i_{pqp^{-1}}(p(h))$. It is easy to check that $H^{G}$ is a crossed group-cograded Hopf quasigroup.
 \\

 Using the mirror reflection technique introduced in Turaev \cite{T}, we can give a construction of crossed group-cograded Hopf quasigroups from a fixed crossed group-cograded Hopf quasigroup as follows.
 \\

 \textbf{Theorem \thesection.5}  Let $(H=\bigoplus_{p\in G}H_{p}, \Delta, \epsilon, S, \pi)$ be a crossed group-cograded Hopf quasigroup, then we can define its mirror $(\widetilde{H}=\bigoplus_{p\in G}\widetilde{H}_{p}, \widetilde{\Delta}, \widetilde{\epsilon}, \widetilde{S}, \widetilde{\pi})$ by the following way:
 \begin{itemize}
   \item[(1)] as an algebra, $\widetilde{H}_{p}=H_{p^{-1}}$, for all $p\in G$;

   \item[(2)] define the comultiplication $\widetilde{\Delta}_{p,q}: \widetilde{H}_{pq}\rightarrow \widetilde{H}_{p}\otimes \widetilde{H}_{q}$ by : for $h_{q^{-1}p^{-1}}\in \widetilde{H}_{pq}$,
   \begin{eqnarray}
   \widetilde{\Delta}_{p, q}(h_{q^{-1}p^{-1}})&=&(\pi_{q}\otimes id_{H_{q^{-1}}})\Delta_{q^{-1}p^{-1}q, q^{-1}}(h_{q^{-1}p^{-1}});\label{th3.52}
   \end{eqnarray}

   \item[(3)] the counit $\widetilde{\epsilon}$ of $\widetilde{H}$ is the original counit $\epsilon$;

   \item[(4)] the antipode $\widetilde{S}_{p}=\pi_{p}S_{p^{-1}}: \widetilde{H}_{p}=H_{p^{-1}}\rightarrow H_{p}=\widetilde{H}_{p^{-1}}$;

   \item[(5)] for all $p\in G$, define the cross action $\widetilde{\pi}_{p}=\pi_{p}$.
 \end{itemize}
 Then $(\widetilde{H}=\bigoplus_{p\in G}\widetilde{H}_{p}, \widetilde{\Delta}, \widetilde{\epsilon}, \widetilde{S}, \widetilde{\pi})$ is also a crossed group-cograded Hopf quasigroup.

 \smallskip\textbf{\emph{Proof}} It is easy to check $\Delta$ is coassociative, and $\epsilon$ is counit of $\widetilde{H}$.

 We will only prove equation (\ref{S1}) of $\widetilde{H}$ holds, the equation (\ref{S2}) of $\widetilde{H}$ is similar. Indeed,
 \begin{eqnarray*}
 &&\mu_{p^{-1}}(id_{\widetilde{H}_{p}}\otimes \mu_{p^{-1}})(\widetilde{S}_{p^{-1}}\otimes id_{\widetilde{H}_{p}}\otimes id_{\widetilde{H}_{p}})(\widetilde{\Delta}_{p^{-1},p}\otimes id_{\widetilde{H}_{p}})\\
 &=&\mu_{p^{-1}}(id_{H_{p^{-1}}}\otimes \mu_{p^{-1}})(\pi_{p^{-1}}S_{p}\otimes id_{H_{p^{-1}}}\otimes id_{H_{p^{-1}}})\big((\pi_{p}\otimes id_{H_{p^{-1}}})\Delta_{p,p^{-1}}\otimes id_{H_{p^{-1}}}\big)\\
 &=&\mu_{p^{-1}}(id_{H_{p^{-1}}}\otimes \mu_{p^{-1}})\big((\pi_{p^{-1}}S_{p}\pi_{p}\otimes id_{H_{p^{-1}}})\Delta_{p,p^{-1}}\otimes id_{H_{p^{-1}}}\big)\\
 &=&\mu_{p^{-1}}(id_{H_{p^{-1}}}\otimes \mu_{p^{-1}})\big((S_{p}\otimes id_{H_{p^{-1}}})\Delta_{p,p^{-1}}\otimes id_{H_{p^{-1}}}\big)\\
 &=&\mu_{p^{-1}}(id_{H_{p^{-1}}}\otimes \mu_{p^{-1}})(S_{p}\otimes id_{H_{p^{-1}}}\otimes id_{H_{p^{-1}}})(\Delta_{p,p^{-1}}\otimes id_{H_{p^{-1}}})\\
 &=&\epsilon \otimes id_{\widetilde{H}_{p}},
 \end{eqnarray*}
 and
 \begin{eqnarray*}
 &&\mu_{p^{-1}}(id_{\widetilde{H}_{p}}\otimes \mu_{p^{-1}})(id_{\widetilde{H}_{p}}\otimes \widetilde{S}_{p^{-1}}\otimes id_{\widetilde{H}_{p}})(\widetilde{\Delta}_{p, p^{-1}}\otimes id_{\widetilde{H}_{p}})\\
 &=&\mu_{p^{-1}}(id_{H_{p^{-1}}}\otimes \mu_{p^{-1}})(id_{H_{p^{-1}}}\otimes \pi_{p^{-1}}S_{p^{-1}}\otimes id_{H_{p^{-1}}})\big((\pi_{p^{-1}}\otimes id_{H_{p}})\Delta_{pp^{-1}p^{-1}, p}\otimes id_{H_{p^{-1}}}\big)\\
 &=&\mu_{p^{-1}}(id_{H_{p^{-1}}}\otimes \mu_{p^{-1}})\big((\pi_{p^{-1}}\otimes S_{p^{-1}pp}\pi_{p^{-1}})\Delta_{p^{-1}, p}\otimes id_{H_{p^{-1}}}\big)\\
 &=&\mu_{p^{-1}}(id_{H_{p^{-1}}}\otimes \mu_{p^{-1}})\big((id_{H_{p^{-1}}}\otimes S_{p})(\pi_{p^{-1}}\otimes \pi_{p^{-1}})\Delta_{p^{-1}, p}\otimes id_{H_{p^{-1}}}\big)\\
 &=&\mu_{p^{-1}}(id_{H_{p^{-1}}}\otimes \mu_{p^{-1}})\big((id_{H_{p^{-1}}}\otimes S_{p})\Delta_{p^{-1}, p}\pi_{p^{-1}}\otimes id_{H_{p^{-1}}}\big)\\
 &=&\mu_{p^{-1}}(id_{H_{p^{-1}}}\otimes \mu_{p^{-1}})(id_{H_{p^{-1}}}\otimes S_{p}\otimes id_{H_{p^{-1}}})(\Delta_{p^{-1}, p}\otimes id_{H_{p^{-1}}})(\pi_{p^{-1}}\otimes id_{H_{p^{-1}}})\\
 &=&(\epsilon\otimes id_{H_{p^{-1}}})(\pi_{p^{-1}}\otimes id_{H_{p^{-1}}})\\
 &=&\epsilon\otimes id_{\widetilde{H}_{p}},
 \end{eqnarray*}
 so the equation (\ref{S1}) of $\widetilde{H}$ holds.

 It is obvious that $\widetilde{\pi}=\pi$ is multiplicative, and each $\pi_{p}$ preserves the counit, so if each $\pi_{p}$ preserves antipode and comultiplication, the mirror $\widetilde{H}$ of $H$ is also a crossed group-cograded Hopf quasigroup. Indeed, for all $p, q\in G$,
 \begin{eqnarray*}
 \widetilde{S}_{pqp^{-1}}\pi_{p}=\pi_{pqp^{-1}}S_{pqp^{-1}}\pi_{p}=\pi_{pqp^{-1}}\pi_{p}S_{q}=\pi_{pq}S_{q}=\pi_{p}\pi_{q}S_{q}=\pi_{p}\widetilde{S}_{q},
 \end{eqnarray*}
 thus $\pi_{p}$ preserves antipode.
 And we finally consider comultiplication, for all $p, q, r \in G$,
 \begin{eqnarray*}
 (\pi_{p}\otimes \pi_{p})\widetilde{\Delta}_{q,r}
 &=&(\pi_{p}\otimes \pi_{p})(\pi_{r}\otimes id_{H_{r^{-1}}})\Delta_{r^{-1}q^{-1}r, r^{-1}}\\
 &=&(\pi_{pr}\otimes \pi_{p})\Delta_{r^{-1}q^{-1}r, r^{-1}},
 \end{eqnarray*}
 and
 \begin{eqnarray*}
 \widetilde{\Delta}_{pqp^{-1}, prp^{-1}}\pi_{p}
 &=&(\pi_{prp^{-1}}\otimes id_{H_{pr^{-1}p^{-1}}})\Delta_{pr^{-1}q^{-1}rp^{-1}, pr^{-1}p^{-1}}\pi_{p}\\
 &=&(\pi_{prp^{-1}}\otimes id_{H_{pr^{-1}p^{-1}}})(\pi_{p}\otimes \pi_{p})\Delta_{r^{-1}q^{-1}r, r^{-1}}\\
 &=&(\pi_{pr}\otimes \pi_{p})\Delta_{r^{-1}q^{-1}r, r^{-1}},
 \end{eqnarray*}
 hence $\pi_{p}$ preserves comultiplication.
 Then we conclude $\widetilde{H}$ is a crossed group-cograded Hopf quasigroup.$\hfill \square$
 \\

 \textbf{Example \thesection.6} Let $H^{G}$ be a crossed group-cograded Hopf quasigroup introduced in Example \thesection.4.
 Set $\widetilde{H}^{G}$ be the same family of algebras $(H_{p}=H)_{p \in G}$ with the same counit,
 the same action $\pi$ of $G$, the comultiplication $\widetilde{\Delta}_{p, q}: H_{pq}\rightarrow H_{p}\otimes H_{q}$,
 and the antipode $\widetilde{S}_{p}: H_{p}\rightarrow H_{p^{-1}}$ defined by
 \begin{eqnarray*}
 &&\widetilde{\Delta}_{p, q}(i_{pq}(h))=\sum_{(h)}i_{p}(q(h_{(e)}))\otimes i_{q}(h_{(2)}),\\
 &&\widetilde{S}_{p}(i_{p}(h))=i_{p^{-1}}(p(S(h)))=i_{p^{-1}}(S(p(h))),
 \end{eqnarray*}
 where $h\in H$. By Theorem \thesection.5 $\widetilde{H}^{G}$ becomes a crossed group-cograded Hopf quasigroup.

 Note that the crossed group-cograded Hopf quasigroups $H^{G}$ and $\widetilde{H}^{G}$, which defined in Example \thesection.3 and \thesection.4 , respectively, are mirrors of each other.
 \\

 \section{ Construction of braided crossed categories}
 \def\theequation{\thesection.\arabic{equation}}
 \setcounter{equation}{0}

 Let $H=\bigoplus_{r\in G}H_{r}$ is a crossed group-cograded Hopf quasigroup with a bijective antipode $S$. We introduce the definition of $p$-Yetter-Drinfeld quasimodules over $H$, then show the category $\mathscr Y\mathscr D\mathscr Q(H)$ of Yetter-Drinfeld quasimodules is a crossed category, and the subcategory $\mathscr Y\mathscr D(H)$ of Yetter-Drinfeld modules over $H$ is a braided crossed category.
 \\

 \textbf{Definition \thesection.1} Let $M$ be a vector space, $(M, \varphi)$ is called a left $H_{p}$-quasimodule if there exists an action $\varphi: H_{p}\otimes M\rightarrow M, h_{p}\otimes m\rightarrow h_{p}\cdot m$ satisfying
 \begin{eqnarray}
 \varphi (\eta_{p}\otimes M)&=&id_{M}, \label{QM1}\\
 \varphi(S_{p^{-1}}\otimes \varphi)(\Delta_{p^{-1}, p}\otimes M)&=&\epsilon\otimes M\nonumber\\
 &=&\varphi (H_{p}\otimes \varphi)(H_{p}\otimes S_{p^{-1}}\otimes M)(\Delta_{p, p^{-1}}\otimes M). \label{QM2}
 \end{eqnarray}

 Using Sweelder notation, for all $h\in H_{e}, m\in M$, (\ref{QM1}) and (\ref{QM2}) is equivalent to
 \begin{eqnarray}
 1_{p}\cdot m&=&m,\label{QM3}\\
 S_{p^{-1}}(h_{(1,p^{-1})})\cdot (h_{(2,p)}\cdot m)&=&\epsilon(h)m \nonumber\\
 &=&h_{(1,p)}\cdot \big(S_{p^{-1}}(h_{2,p^{-1}})\cdot m\big).\label{QM4}
 \end{eqnarray}

 Moreover, if the condition (\ref{QM2}) is instead by $h\cdot(g\cdot m)=(hg)\cdot m$, where $h, g\in H_{p}$, then left $H_{p}$-quasimodule is a left $H_{p}$-module.
 \\

 \textbf{Definition \thesection.2} Let $V$  be a vector space and $p$ a fixed element in group $G$. A left-right $p$-Yetter-Drinfeld quasimodule is a couple $(V, \rho^{V}=(\rho_{r}^{V})_{r\in G})$, where $V$ is a unital $H_{p}$-quasimodule, and for any $r\in G, \rho_{r}^{V}: V\rightarrow V\otimes H_{r}$ is a $k$-linear morphism, denoted by $\rho_{r}^{V}(v)=v_{(0)}\otimes v_{(1,r)}$ such that the following conditions are satisfied:
 \begin{itemize}
   \item[(1)] V is coassociative in the sense that, for any $r_{1},r_{2}\in G$, we have
   \begin{eqnarray*}
   (\rho_{r_{1}}^{V}\otimes id_{r_{2}})\rho_{r_{2}}^{V}=(id_{V}\otimes \Delta_{r_{1},r_{2}})\rho_{r_{1},r_{2}}^{V};
   \end{eqnarray*}
   \item[(2)] V is counitary in the sense that
   \begin{eqnarray*}
   (id_{V}\otimes \epsilon)\rho_{e}^{V}=id_{V};
   \end{eqnarray*}
   \item[(3)] V is crossed in the sense that for all $v\in V$, $r\in G$ and $h,g\in H_{(r)}$,
   \begin{eqnarray}
   &h_{(1,p)}\cdot v_{(0)}\otimes h_{(2,r)}v_{(1,r)}=(h_{(2,p)}\cdot v)_{(0)}\otimes (h_{(2,p)}\cdot v)_{(1,r)}\pi_{p^{-1}}(h_{(1,prp^{-1})}), \label{JR1}\\
   &v_{(0)}\otimes v_{(1,r)}(hg)=v_{(0)}\otimes (v_{(1,r)}h)g, \label{JH1}\\
   &v_{(0)}\otimes (hv_{(1,r)})g=v_{(0)}\otimes h(v_{(1,r)}g). \label{JH2}
   \end{eqnarray}
 \end{itemize}

 Given two $p$-Yetter-Drinfeld quasimodules $(V, \rho{}^{V})$ and $(W, \rho^{W})$, a morphism of this two $p$-Yetter-Drinfeld quasimodules $f: (V, \rho^{V})\rightarrow (W,\rho^{W})$ is an $H_{p}$-linear map $f: V\rightarrow W$ such that for any $r\in G$,
 \begin{eqnarray*}
 \rho_{r}^{W}f=(f\otimes id_{r})\rho_{r}^{V},
 \end{eqnarray*}
 i.e., for all $v\in V$,
 \begin{eqnarray*}
 f(v)_{(0)}\otimes f(v)_{(1,r)}=f(v_{(0)})\otimes v_{(1,r)}.
 \end{eqnarray*}

 Then we have the category $\mathscr Y\mathscr D\mathscr Q(H)_{p}$ of $p$-Yetter-Drinfeld quasimodules, the composition of morphisms of $p$-Yetter-Drinfeld quasimodules is the standard composition of the underlying linear maps. Moreover, if we assume that $V$ is a left $H_{p}$-module, then we say that is a left-right $p$-Yetter-Drinfeld module. Obviously, left-right $p$-Yetter-Drinfeld modules with the obvious morphisms is a subcategory of $\mathscr Y\mathscr D\mathscr Q(H)_{p}$, and denote it by $\mathscr Y\mathscr D(H)_{p}$.
 \\

 \textbf{Proposition \thesection.3} The equation (\ref{JR1}) is equivalent to
 \begin{eqnarray}
 (h_{p}\cdot v)_{(0)}\otimes (h_{p}\cdot v)_{(1,r)}=h_{(2,p)}\cdot v_{(0)}\otimes \big(h_{(3,r)}v_{(1,r)}\big)S^{-1}\pi_{p^{-1}}(h_{(1,pr^{-1}p^{-1})}), \label{JR2}
 \end{eqnarray}
 for all $h_{p}\in H_{p}$ and $v\in V$.

 \textbf{\emph{Proof}}
 Indeed, suppose the condition (\ref{JR1}) holds, then we have
 \begin{eqnarray*}
 &&h_{(2,p)}\cdot v_{(0)}\otimes \big(h_{(3,r)}v_{(1,r)}\big)S^{-1}\pi_{p^{-1}}(h_{(1,pr^{-1}p^{-1})})\\
 &=&(h_{(3,p)}\cdot v)_{(0)}\otimes \big((h_{(3,p)}\cdot v)_{(1,r)}\pi_{p^{-1}}(h_{(2,prp^{-1})})\big)S^{-1}\pi_{p^{-1}}(h_{(1,pr^{-1}p^{-1})})\\
 &=&(h_{(3,p)}\cdot v)_{(0)}\otimes (h_{(3,p)}\cdot v)_{(1,r)}\Big(\pi_{p^{-1}}\big(h_{(2,prp^{-1})}S^{-1}(h_{(1,pr^{-1}p^{-1})})\big)\Big)\\
 &=&(h_{(2,p)}\cdot v)_{(0)}\otimes (h_{(2,p)}\cdot v)_{(1,r)}\pi_{p^{-1}}\epsilon(h_{e})\\
 &=&(h_{(2,p)}\cdot v)_{(0)}\otimes (h_{(2,p)}\cdot v)_{(1,r)}\epsilon(h_{e})\\
 &=&(h_{p}\cdot v)_{(0)}\otimes (h_{p}\cdot v)_{(1,r)}
 \end{eqnarray*}
 where the first equality follows by (\ref{JR1}), the others relies on the properties of crossed group-cograded Hopf quasigroup.

 Conversely, if the equation (\ref{JR2}) holds, then
 \begin{eqnarray*}
 &&(h_{(2,p)}\cdot v)_{(0)}\otimes (h_{(2,p)}\cdot v)_{(1,r)}\pi_{p^{-1}}(h_{(1,prp^{-1})})  \\
 &=&h_{(3,p)}\cdot v_{(0)}\otimes \big(h_{(4,r)}v_{(1,r)}\big)S^{-1}\pi_{p^{-1}}(h_{(2,prp^{-1})})\pi_{p^{-1}}(h_{(1,prp^{-1})})\\
 &=&h_{(3,p)}\cdot v_{(0)}\otimes \big(h_{(4,r)}v_{(1,r)}\big)\pi_{p^{-1}}(S^{-1}(h_{(2,prp^{-1})})h_{(1,prp^{-1})})\\
 &=&h_{(1,p)}\cdot v_{(0)}\otimes \big(h_{(2,r)}v_{(1,r)}\big)\pi_{p^{-1}}\epsilon(h_{e})\\
 &=&h_{(1,p)}\cdot v_{(0)}\otimes h_{(2,r)}v_{(1,r)}
 \end{eqnarray*}
 where the first equality follows by (\ref{JR2}), the rest are follows by the properties of crossed group-cograded Hopf quasigroup. $\hfill \square$
 \\

 \textbf{Remark} According to the equation (\ref{JH2}), the condition (\ref{JR2}) is equivalent to
 \begin{eqnarray}
 (h_{p}\cdot v)_{(0)}\otimes (h_{p}\cdot v)_{(1,r)}=h_{(2,p)}\cdot v_{(0)}\otimes h_{(3,r)}\big(v_{(1,r)}S^{-1}\pi_{p^{-1}}(h_{(1,pr^{-1}p^{-1})})\big).\label{JR3}
 \end{eqnarray}

 \textbf{Proposition \thesection.4} If $(V, \rho^{V})\in \mathscr Y\mathscr D\mathscr Q(H)_{p}$ and $(W, \rho^{W})\in \mathscr Y\mathscr D\mathscr Q(H)_{q}$, then $V\otimes W\in \mathscr Y\mathscr D\mathscr Q(H)_{pq}$ with the module and comodule structures as follows:
 \begin{eqnarray}
 h_{pq}\cdot (v\otimes w)&=&h_{(1,p)}\cdot v\otimes h_{(2,q)}\cdot w, \label{HHQ1}\\
 \rho_{r}^{V\otimes W}(v\otimes w)&=&v_{(0)}\otimes w_{(0)}\otimes w_{(1,r)}\pi_{q^{-1}}(v_{(1,qrq^{-1})}) \label{HHQ2},
 \end{eqnarray}
 where $v\in V, w\in W$ and $h_{pq}\in H_{pq}$.

 \textbf{\emph{Proof}} We first check that $V\otimes W$ is a left $H_{pq}$-quasimodule, and the unital property is obvious.
 We only check the left hand side of equation (\ref{QM2}), the right hand is similar. For all $v\in V, w\in W$,
 \begin{eqnarray*}
 &&h_{(1,pq)}\cdot \big(S^{-1}(h_{(2,(pq)^{-1})})\cdot (v\otimes w)\big)\\
 &=&h_{(1,pq)}\cdot \big(S^{-1}(h_{(2, p^{-1})})\cdot v \otimes S^{-1}(h_{(2,q^{-1})})\cdot w\big)\\
 &=&\Big(h_{(1,p)}\cdot \big(S^{-1}(h_{(2, p^{-1})})\cdot v\big)\Big)\otimes \Big(h_{(2,q)}\cdot \big(S^{-1}(h_{(2,q^{-1})})\cdot w\big)\Big)\\
 &=&\big(\epsilon (h_{(1,e)})\cdot v\big)\otimes \big(\epsilon (h_{(2,e)})\cdot w\big)\\
 &=&\epsilon (h_{e})\cdot (v\otimes w)
 \end{eqnarray*}
 where the first and second equalities are rely on (\ref{HHQ1}), the third equality is follows by (\ref{QM2}). Since then $V\otimes W$ is a left $H_{pq}$-quasimodule.

 In the following equations ,we check that coassociative condition holds:
 \begin{eqnarray*}
 &&(id_{V\otimes W}\otimes \Delta_{r_{1}r_{2}})\rho_{r_{1}r_{2}}(v\otimes w)\\
 &=&(id_{V\otimes W}\otimes \Delta_{r_{1}r_{2}})\big(v_{(0)}\otimes w_{(0)}\otimes w_{(1,r_{1}r_{2})}\pi_{q^{-1}}(v_{(1,qr_{1}r_{2}q^{-1})})\big)\\
 &=&v_{(0)}\otimes w_{(0)}\otimes w_{(1,r_{1})}\pi_{q^{-1}}v_{(1,qr_{1}q^{-1})}\otimes w_{(2,r_{2})}\pi_{q^{-1}}(v_{(2,qr_{2}q^{-1})}),
 \end{eqnarray*}
 and
 \begin{eqnarray*}
 &&(\rho_{r_{1}}\otimes id_{r_{2}})\rho_{r_{2}}(v\otimes w)\\
 &=&(\rho_{r_{1}}\otimes id_{r_{2}})(v_{(0)}\otimes w_{(0)}\otimes w_{(1, r_{2})}\pi_{q^{-1}}(v_{1,qr_{2}q^{-1}}))\\
 &=&v_{(0)(0)}\otimes w_{(0)(0)}\otimes w_{(0)(1, r_{1})}\pi_{q^{-1}}(v_{(0)(1,qr_{1}q^{-1})})\otimes w_{(1,r_{2})}\pi_{q^{-1}}(v_{1,qr_{2}q^{-1}})\\
 &=&v_{(0)}\otimes w_{(0)}\otimes w_{(1,r_{1})}\pi_{q^{-1}}(v_{(1,qr_{1}q^{-1})})\otimes w_{(2,r_{2})}\pi_{q^{-1}}(v_{(2,qr_{2}q^{-1})}).
 \end{eqnarray*}
 This shows that $(id_{V\otimes W}\otimes \Delta_{r_{1}r_{2}})\rho_{r_{1}r_{2}}=(\rho_{r_{1}}\otimes id_{r_{2}})\rho_{r_{2}}$.

 And the counitary condition is easy to show. Then we check the crossed condition as follows:
 \begin{eqnarray*}
 &&h_{(1,pq)}\cdot (v\otimes w)_{(0)}\otimes h_{(2,r)}(v\otimes w)_{(1,r)}\\
 &\stackrel{(\ref{HHQ2})}=&h_{(1,pq)}\cdot (v_{(0)}\otimes w_{(0)})\otimes h_{(2,r)}\big(w_{(1,r)}\pi_{q^{-1}}(v_{(1,qrq^{-1})})\big)\\
 &\stackrel{(\ref{JH1})}=&h_{(1,pq)}\cdot (v_{(0)}\otimes w_{(0)})\otimes \big(h_{(2,r)}w_{(1,r)}\big)\pi_{q^{-1}}(v_{(1,qrq^{-1})})\\
 &\stackrel{(\ref{HHQ1})}=&h_{(1,p)}\cdot v_{(0)}\otimes h_{(2,q)}\cdot w_{(0)}\otimes \big(h_{(3,r)}w_{(1,r)}\big)\pi_{q^{-1}}(v_{1,qrq^{-1}})\\
 &\stackrel{(\ref{JR1})}=&h_{(1,p)}\cdot v_{(0)}\otimes (h_{(3,q)}\cdot w)_{(0)}\otimes (h_{(3,q)}\cdot w)_{(1,r)}\pi_{q^{-1}}(h_{(2,qrq^{-1})})\pi_{q^{-1}}(v_{(1,qrq^{-1})})\\
 &=&h_{(1,p)}\cdot v_{(0)}\otimes (h_{(3,q)}\cdot w)_{(0)}\otimes (h_{(3,q)}\cdot w)_{(1,r)}\pi_{q^{-1}}(h_{(2,qrq^{-1})}v_{(1,qrq^{-1})})\\
 &\stackrel{(\ref{JR1})}=&(h_{(2,p)}\cdot v)_{(0)}\otimes (h_{(3,q)}\cdot w)_{(0)}\\
 &&\otimes (h_{(3,q)}\cdot w)_{(1,r)}\pi_{q^{-1}}\big((h_{(2,q)}\cdot v)_{(1,qrq^{-1})}\pi_{p^{-1}}(h_{(1,pqrq^{-1}p^{-1})})\big)\\
 &=&(h_{(2,p)}\cdot v)_{(0)}\otimes (h_{(3,q)}\cdot w)_{(0)}\\
 &&\otimes (h_{(3,q)}\cdot w)_{(1,r)}\pi_{q^{-1}}\big((h_{(2,q)}\cdot v)_{(1,qrq^{-1})}\big)\pi_{q^{-1}p^{-1}}(h_{(1,pqrq^{-1}p^{-1})})\\
 &\stackrel{(\ref{HHQ2})}=&(h_{(2,p)}\cdot v \otimes h_{(3,q)}\cdot w)_{(0)} \otimes (h_{(2,p)}\cdot v \otimes h_{(3,q)}\cdot w)_{(1,r)}\pi_{q^{-1}p^{-1}}(h_{(1,pqrq^{-1}p^{-1})})\\
 &\stackrel{(\ref{HHQ1})}=&h_{(2,pq)}\cdot (v\otimes w)_{(0)}\otimes \big(h_{(2,pq)}\cdot (v\otimes w)\big)_{(1,r)}\pi_{q^{-1}p^{-1}}(h_{(1,pqrq^{-1}p^{-1})}).
 \end{eqnarray*}

 Finally, we check the equation (\ref{JH1}), and the equation (\ref{JH2}) is similar.
 \begin{eqnarray*}
 (v\otimes w)_{(0)}\otimes (v\otimes w)_{(1,r)}(hg)
 &=&v_{(0)}\otimes w_{(0)}\otimes \big(w_{(1,r)}\pi_{q^{-1}}(v_{(1,qrq^{-1})})\big)(hg)\\
 &=&v_{(0)}\otimes w_{(0)}\otimes w_{(1,r)}\big(\pi_{q^{-1}}(v_{(1,qrq^{-1})})(hg)\big)\\
 &=&v_{(0)}\otimes w_{(0)}\otimes w_{(1,r)}\Big(\big(\pi_{q^{-1}}(v_{(1,qrq^{-1})})h\big)g\Big)\\
 &=&v_{(0)}\otimes w_{(0)}\otimes \Big(w_{(1,r)}\big(\pi_{q^{-1}}(v_{(1,qrq^{-1})})h\big)\Big)g\\
 &=&v_{(0)}\otimes w_{(0)}\otimes \Big(\big(w_{(1,r)}\pi_{q^{-1}}(v_{(1,qrq^{-1})})\big)h\Big)g\\
 &=&(v\otimes w)_{(0)}\otimes \big((v\otimes w)_{(1,r)}h\big)g.
 \end{eqnarray*}
 Hence $V\otimes W\in \mathscr Y\mathscr D\mathscr Q(H)_{pq}$. $\hfill \square$
 \\

 Following Turaev's left index notation, let $V\in \mathscr Y \mathscr D\mathscr Q(H)_{p}$, the object $^{q}V$ has the same underlying vector space as V. Given $v\in V$, we denote $^{q}v$ the corresponding element in $^{q}V$.

 \smallskip\textbf{Proposition \thesection.5} Let $(V, \rho^{V})\in \mathscr Y \mathscr D \mathscr Q(H)_{p}$ and $q\in G$. Set $^{q}V=V$ as a vector space with structures
 \begin{eqnarray}
 &&h_{qpq^{-1}}\cdot {}^{q}v={}^{q}(\pi_{q^{-1}}(h_{qpq^{-1}})\cdot v) \label{qa}\\
 &&\rho_{r}^{^{q}V}(^{q}v)={}^{q}(v_{(0)})\otimes \pi_{q}(v_{(1,q^{-1}rq)}) \label{qco}
 \end{eqnarray}
 for any $v\in V$ and $h_{qpq^{-1}}\in H_{qpq^{-1}}$. Then ${}^{q}V\in \mathscr Y \mathscr D \mathscr Q(H)_{qpq^{-1}}$.

 \textbf{\emph{Proof}} We first cheek that $^{q}V$ is a left $H_{qpq^{-1}}$-quasimodule. The condition (\ref{QM1}) is easy to check.  Next,  we prove the condition (\ref{QM2}).
 \begin{eqnarray*}
 h_{(1,qpq^{-1})}\cdot (S^{-1}(h_{(2,qp^{-1}q^{-1})})\cdot {}^{q}v)
 &=&h_{(1,qpq^{-1})}\cdot \Big({}^{q}\big(\pi_{q^{-1}}\big(S^{-1}(h_{(2,qp^{-1}q^{-1})})\big)\cdot v\big)\Big)\\
 &=&{}^{q}\Big(\pi_{q^{-1}}(h_{(1,qpq^{-1})})\cdot \big(\pi_{q^{-1}}\big(S^{-1}(h_{(2,qp^{-1}q^{-1})})\big)\cdot v\big)\Big)\\
 &=&{}^{q}\Big(\epsilon\big(\pi_{q^{-1}}(h_{e})\big)\cdot v\Big)\\
 &=&\epsilon\big(\pi_{q^{-1}}(h_{e})\big)\cdot {}^{q}v\\
 &=&\epsilon(h_{e}){}^{q}v.
 \end{eqnarray*}
 The proof of the other hand is similar as above, so $^{q}V$ is a left $H_{qpq^{-1}}$-quasimodule. And the coassociative and counitary are also satisfied.

 In the following, we show that the crossing condition holds:
 \begin{eqnarray*}
 &&(h_{qpq^{-1}}\cdot {}^{q}v)_{(0)}\otimes (h_{qpq^{-1}}\cdot {}^{q}v)_{(1,r)}\\
 &\stackrel{(\ref{qa})}=&\Big({}^{q}\big(\pi_{q^{-1}}(h_{qpq^{-1}})\cdot v\big)\Big)_{(0)}\otimes \Big({}^{q}\big(\pi_{q^{-1}}(h_{qpq^{-1}})\cdot v\big)\Big)_{(1,r)}\\
 &\stackrel{(\ref{qco})}=&{}^{q}\Big(\big(\pi_{q^{-1}}(h_{qpq^{-1}})\cdot v\big)_{(0)}\Big)\otimes \pi_{q}\Big(\big(\pi_{q^{-1}}(h_{qpq^{-1}})\big)_{(1,r)}\Big)\\
 &\stackrel{(\ref{JR2})}=&{}^{q}\Big(\pi_{q^{-1}}(h_{qpq^{-1}})_{(2,p)}\cdot v_{(0)}\Big)\\
 &&\otimes \pi_{q}\Big(\big(\pi_{q^{-1}}(h_{qpq^{-1}})_{(3,q^{-1}rq)}v_{(1,q^{-1}rq)}\big)S^{-1}\pi_{p^{-1}}\big(\pi_{q^{-1}}(h_{qpq^{-1}})_{(1,pq^{-1}r^{-1}qp^{-1})}\big)\Big)\\
 &=&{}^{q}\Big(\pi_{q^{-1}}(h_{(2,qpq^{-1})})\cdot v_{(0)}\Big)\\
 &&\otimes \pi_{q}\Big(\big(\pi_{q^{-1}}(h_{(3,r)})v_{(1,q^{-1}rq)}\big)S^{-1}\pi_{p^{-1}}\big(\pi_{q^{-1}}(h_{(1,qpq^{-1}r^{-1}qp^{-1}q^{-1})})\big)\Big)\\
 &=&{}^{q}\Big(\pi_{q^{-1}}(h_{(2,qpq^{-1})})\cdot v_{(0)}\Big)\otimes \big(h_{(3,r)}\pi_{q}(v_{(1,q^{-1}rq)})\big)S^{-1}\pi_{qpq^{-1}}(h_{(1,qpq^{-1}r^{-1}qp^{-1}q^{-1})})\\
 &\stackrel{(\ref{qa})}=&h_{(2,qpq^{-1})}\cdot {}^{q}(v_{(0)})\otimes  \big(h_{(3,r)}\pi_{q}(v_{(1,q^{-1}rq)})\big)S^{-1}\pi_{qpq^{-1}}(h_{(1,qpq^{-1}r^{-1}qp^{-1}q^{-1})})\\
 &\stackrel{(\ref{qco})}=&h_{(2,qpq^{-1})}\cdot({}^{q}v)_{(0)}\otimes \big(h_{(3,r)}({}^{q}v)_{(1,r)}\big)S^{-1}\pi_{qpq^{-1}}(h_{(1,qpq^{-1}r^{-1}qp^{-1}q^{-1})}).
 \end{eqnarray*}

 At last, we will check that the quasimodule coassociative conditions are hold. We just compute the equation (\ref{JH1}), the equation (\ref{JH2}) is similar. For all ${}^{q}v \in {}^{q}V,\, h,g \in H_{r}$,
 \begin{eqnarray*}
 ({}^{q}v)_{(0)}\otimes ({}^{q}v)_{(1,r)}(hg)
 &=&{}^{q}(v_{(0)})\otimes \pi_{q}(v_{(1,q^{-1}rq)})(hg)\\
 &=&{}^{q}(v_{(0)})\otimes \big(\pi_{q}(v_{(1,q^{-1}rq)})h\big)g\\
 &=&({}^{q}v)_{(0)}\otimes \big(({}^{q}v)_{(1,r)}h\big)g,
 \end{eqnarray*}
 where the first and third equalities are rely on (\ref{qco}), the second one is follows by (\ref{JH1}).
 This completes the proof.$\hfill \square$
 \\

 \textbf{Proposition \thesection.6} Let $(V,\rho^{V})\in \mathscr Y\mathscr D\mathscr Q(H)_{p}$ and $(W,\rho^{W})\in \mathscr Y\mathscr D\mathscr Q(H)_{q}$. Then $^{st}V={}^{s}(^{t}V)$ as an object in $\mathscr Y\mathscr D\mathscr Q(H)_{stpt^{-1}s^{-1}}$, and $^{s}(V\otimes W)={}^{s}V\otimes {}^{s}W$ as an object in $\mathscr Y\mathscr D\mathscr Q(H)_{spqs^{-1}}$.

 \textbf{\emph{Proof}} We first check that $^{st}V={}^{s}(^{t}V)$ as an object in $\mathscr Y\mathscr D\mathscr Q(H)_{stpt^{-1}s^{-1}}$. It is obviously that both $^{st}V$ and ${}^{s}(^{t}V)$ are in the category $\mathscr Y\mathscr D\mathscr Q(H)_{stpt^{-1}s^{-1}}$. Then we show that the action and coaction of these two $stpt^{-1}s^{-1}$-Yetter-Drinfeld quasimodules are exactly equivalent.

 As $^{st}V$ is a $stpt^{-1}s^{-1}$-Yetter-Drinfeld quasimodules with the structures
 \begin{eqnarray*}
 h_{stpt^{-1}s^{-1}}\cdot {}^{st}v&=&{}^{st}\big(\pi_{t^{-1}s^{-1}}(h_{stpt^{-1}s^{-1}})\cdot v\big),\\
 \rho^{{}^{st}V}_{r}({}^{st}v)&=&{}^{st}(v_{(0)})\otimes \pi_{st}(v_{(1,t^{-1}s^{-1}rst)}).
 \end{eqnarray*}
 Then we show ${}^{s}(^{t}V)$ is $stpt^{-1}s^{-1}$-Yetter-Drinfeld quasimodules with the same structures of $^{st}V$.
 Indeed, the action of ${}^{s}(^{t}V)$ is
 \begin{eqnarray*}
 h_{stpt^{-1}s^{-1}}\cdot {}^{s}({}^{t}v)
 &=&{}^{s}\big(\pi_{s^{-1}}(h_{stpt^{-1}s^{-1}})\cdot {}^{t}v\big)\\
 &=&{}^{s}\Big({}^{t}\big(\pi_{t^{-1}}\pi_{s^{-1}}(h_{stpt^{-1}s^{-1}})\cdot v\big)\Big)\\
 &=&{}^{s}{}^{t}\big(\pi_{t^{-1}s^{-1}}(h_{stpt^{-1}s^{-1}})\cdot v\big).
 \end{eqnarray*}
 Hence ${}^{s}(^{t}V)$ has the same cation with $^{st}V$.

 And the coaction of ${}^{s}(^{t}V)$ is
 \begin{eqnarray*}
 \rho^{{}^{s}(^{t}V)}\big({}^{s}(^{t}v)\big)
 &=&{}^{s}\big(({}^{t}v)_{(0)}\big)\otimes \pi_{s}\big(({}^{t}v)_{(1, s^{-1}rs)}\big)\\
 &=&{}^{s}\big({}^{t}(v_{(0)})\big)\otimes \pi_{s}\big(\pi_{t}(v_{(1, t^{-1}s^{-1}rst)})\big)\\
 &=&{}^{st}(v_{(0)})\otimes \pi_{st}(v_{(1, t^{-1}s^{-1}rst)}).
 \end{eqnarray*}
 Hence $^{st}V={}^{s}(^{t}V)$ as an object in $\mathscr Y\mathscr D\mathscr Q(H)_{stpt^{-1}s^{-1}}$.

 As $^{s}(V\otimes W)$ is a $spqs^{-1}$-Yetter-Drinfeld quasimodules with the structures
 \begin{eqnarray*}
 h_{spqs^{-1}}\cdot {}^{s}(v\otimes w)&=&{}^{s}\big(\pi_{s^{-1}}(h_{spqs^{-1}})\cdot (v\otimes w)\big)\\
 &=&{}^{s}(\pi_{s^{-1}}(h_{(1,sps^{-1})})\cdot v\otimes \pi_{s^{-1}}(h_{(2,sqs^{-1})})\cdot w),\\
 \rho^{{}^{s}(V\otimes W)}_{r}\big({}^{s}(v\otimes w)\big)&=&{}^{s}\big((v\otimes w)_{(0)}\big)\otimes \pi_{s}\big((v\otimes w)_{(1,s^{-1}rs)}\big)\\
 &=&{}^{s}(v_{(0)}\otimes w_{(0)})\otimes \pi_{s}\big(w_{(1,s^{-1}rs)}\pi_{q^{-1}}(v_{(1,qs^{-1}rsq^{-1})})\big)\\
 &=&{}^{s}(v_{(0)}\otimes w_{(0)})\otimes \pi_{s}(w_{(1,s^{-1}rs)})\pi_{sq^{-1}}(v_{(1,qs^{-1}rsq^{-1})}).
 \end{eqnarray*}

 Then we show ${}^{s}V\otimes {}^{s}W$ is a $spqs^{-1}$-Yetter-Drinfeld quasimodules with the same structures of $^{s}(V\otimes W)$.
 Indeed, the action of ${}^{s}V\otimes {}^{s}W$ is
 \begin{eqnarray*}
 h_{spqs^{-1}}\cdot ({}^{s}v\otimes {}^{s}w)
 &=&h_{(1,sps^{-1})}\cdot {}^{s}v\otimes h_{(2,sqs^{-1})}\cdot {}^{s}w\\
 &=&{}^{s}(\pi_{s^{-1}}(h_{(1,sps^{-1})})\cdot v)\otimes {}^{s}(\pi_{s^{-1}}(h_{(2,sqs^{-1})})\cdot w)\\
 &=&{}^{s}(\pi_{s^{-1}}(h_{(1,sps^{-1})})\cdot v\otimes \pi_{s^{-1}}(h_{(2,sqs^{-1})})\cdot w)
 \end{eqnarray*}
 Hence ${}^{s}V\otimes {}^{s}W$ has the same cation with $^{s}(V\otimes W)$.

 And the coaction of ${}^{s}V\otimes {}^{s}W$ is
 \begin{eqnarray*}
 \rho^{{}^{s}V\otimes {}^{s}W}_{r}({}^{s}v\otimes {}^{s}w)
 &=&({}^{s}v)_{(0)}\otimes ({}^{s}w)_{(0)}\otimes ({}^{s}w)_{(1,r)}\pi_{sq^{-1}s^{-1}}\big(({}^{s}v)_{(1,sqs^{-1}rsq^{-1}s^{-1})}\big)\\
 &=&({}^{s}v)_{(0)}\otimes {}^{s}(w_{(0)})\otimes \pi_{s}(w_{(1,s^{-1}rs)})\pi_{sq^{-1}s^{-1}}\big(\pi_{s}(v_{(1,qs^{-1}rsq^{-1})})\big)\\
 &=&{}^{s}(v_{(0)})\otimes {}^{s}(w_{(0)})\otimes \pi_{s}(w_{(1,s^{-1}rs)})\pi_{sq^{-1}}(v_{(1,qs^{-1}rsq^{-1})})\\
 &=&{}^{s}(v_{(0)}\otimes w_{(0)})\otimes \pi_{s}(w_{(1,s^{-1}rs)})\pi_{sq^{-1}}(v_{(1,qs^{-1}rsq^{-1})}).
 \end{eqnarray*}
 Thus $^{s}(V\otimes W)={}^{s}V\otimes {}^{s}W$ as an object in $\mathscr Y\mathscr D\mathscr Q(H)_{spqs^{-1}}$.$\hfill \square$
 \\

 For a crossed group-cograded Hopf quasigroup $H$, we define $\mathscr Y\mathscr D\mathscr Q(H)$ as the disjoint union of all $\mathscr Y\mathscr D\mathscr Q(H)_{p}$ with $p\in G$.
 If we endow $\mathscr Y\mathscr D\mathscr Q(H)$ with tensor product as in Proposition \thesection.4, then we get the the following result.
 \\

 \textbf{Theorem \thesection.7} The Yetter-Drinfeld quasimodules category  $\mathscr Y\mathscr D\mathscr Q(H)$ is a crossed category.

 \textbf{\emph{Proof}} By Proposition \thesection.5 we can give a group homomorphism $\phi: G\rightarrow Aut(\mathscr Y\mathscr D\mathscr Q(H))$, $p\longmapsto \phi_{p}$ by
 \begin{eqnarray*}
 \phi_{p}: \mathscr Y\mathscr D\mathscr Q(H)_{q}\rightarrow \mathscr Y\mathscr D\mathscr Q(H)_{pqp^{-1}}, \qquad \phi_{p}(W)={}^{p}W,
 \end{eqnarray*}
 where the functor $\phi_{p}$ act as follows: given a morphism $f:(V,\rho^{V})\rightarrow (W,\rho^{W})$, for any $v\in V$, we set $(^{p}f)(^{p}v)={}^{p}(f(v))$.

 Then it is easy to prove $\mathscr Y\mathscr D\mathscr Q(H)$ is a crossed category. $\hfill \square$
 \\

 Following the ideas by $\rm{\acute{A}}$lonso in \cite{AFG15}, we will consider $\mathscr Y\mathscr D(H)_{p}$ the category of left-right $p$-Yetter-Drinfeld modules over $H$, which is a subcategory of $\mathscr Y\mathscr D\mathscr Q(H)_{p}$.
 \\

 \textbf{Proposition \thesection.8} Let $(V,\rho^{V})\in \mathscr Y\mathscr D(H)_{p}$ and $(W,\rho^{W})\in \mathscr Y\mathscr D(H)_{q}$. Set ${}^{V}W={}^{p}W$ as an object in $\mathscr Y\mathscr D(H)_{pqp^{-1}}$. Define the map
 \begin{eqnarray}
 &&C_{V,W}: V\otimes W\rightarrow {}^{V}W\otimes V \nonumber \\
 &&C_{V,W}(v\otimes w)={}^{p}\big(S_{q^{-1}}(v_{(1,q^{-1})})\cdot w\big)\otimes v_{(0)} \label{BZ}
 \end{eqnarray}
 Then $C_{V,W}$ is $H$-linear, $H$-colinear and satisfies the conditions:
 \begin{eqnarray}
 &&C_{V\otimes W, X}=(C_{V,^{W}X}\otimes id_{W})(id_{V}\otimes C_{W, X})\\
 &&C_{V, W\otimes X}=(id_{^{V}W}\otimes C_{V,X})(C_{V,W}\otimes id_{X})
 \end{eqnarray}
 for $X\in \mathscr Y\mathscr D(H)_{s}$. Moreover, $C_{^{s}V,^{s}W}={}^{s}(\cdot) C_{V,W}$.

 \textbf{\emph{Proof}} We first show that $C_{V,W}$ is $H$-linear. First, compute
 \begin{eqnarray*}
 &&C_{V,W}\big(h_{pq}\cdot (v\otimes w)\big)\\
 &\stackrel{(\ref{HHQ1})}=&C_{V,W}(h_{(1,p)}\cdot v\otimes h_{(2,q)}\cdot w)\\
 &\stackrel{(\ref{BZ})}=&^{p}\Big(S_{q^{-1}}\big((h_{(1,p)}\cdot v)_{(1,q^{-1})}\big)\cdot (h_{(2,q)}\cdot w)\Big)\otimes (h_{(1,p)}\cdot v)_{(0)}\\
 &\stackrel{(\ref{JR2})}=&^{p}\Big(S_{q^{-1}}\big(h_{(3,q^{-1})}v_{(1,q^{-1})}S^{-1}\pi_{p^{-1}}(h_{(1,pqp^{-1})})\big)\cdot (h_{(4,q)}\cdot w)\Big)\otimes \big(h_{(2,p)}\cdot v_{(0)}\big)\\
 &=&^{p}\Big(\big(\pi_{p^{-1}}(h_{(1,pqp^{-1})})S_{q^{-1}}(v_{(1,q^{-1})})S^{-1}(h_{3,q^{-1}})\big)\cdot (h_{(4,q)}\cdot w)\Big)\otimes \big(h_{(2,p)}\cdot v_{(0)}\big)\\
 &=&^{p}\Big(\big(\pi_{p^{-1}}(h_{(1,pqp^{-1})})S_{q^{-1}}(v_{(1,q^{-1})})\big)\cdot \big(S^{-1}(h_{(3,q^{-1})})h_{(4,q)}\cdot w\big)\Big)\otimes \big(h_{(2,p)}\cdot v_{(0)}\big)\\
 &=&^{p}\Big(\pi_{p^{-1}}(h_{(1,pqp^{-1})})\big(S_{q^{-1}}(v_{(1,q^{-1})})\cdot w\big)\Big)\otimes \big(h_{(2,p)}\cdot v_{(0)}\big)\\
 &\stackrel{(\ref{qa})}=&\Big(h_{(1,pqp^{-1})}\cdot {}^{p}\big(S_{q^{-1}}(v_{(1,q^{-1})})\cdot w\big)\Big)\otimes \big(h_{(2,p)}\cdot v_{(0)}\big)\\
 &\stackrel{(\ref{HHQ1})}=&h_{pq}\cdot \Big(^{p}\big(S_{q^{-1}}(v_{(1,q^{-1})})\cdot w\big)\otimes v_{(0)}\Big)\\
 &\stackrel{(\ref{BZ})}=&h_{pq}\cdot C_{V,W}(v\otimes w),
 \end{eqnarray*}
 so we have $C_{V,W}\big(h_{pq}\cdot (v\otimes w)\big)=h_{pq}\cdot C_{V,W}(v\otimes w)$, that is, $C_{V,W}$ is $H$-linear.

 Secondly, we prove that $C_{V,W}$ is $H$-colinear. In fact,
 \begin{eqnarray*}
 &&\rho_{r}^{^{V}W\otimes V}C_{V,W}(v\otimes w)\\
 &=&\rho_{r}^{^{V}W\otimes V}\Big(^{p}\big(S_{q^{-1}}(v_{(1,q^{-1})})\cdot w\big)\otimes v_{(0)}\Big)\\
 &\stackrel{(\ref{HHQ2})}=&^{p}\big(S_{q^{-1}}(v_{(1,q^{-1})})\cdot w\big)_{(0)}\otimes v_{(0)(0)}\otimes v_{(0)(1,r)}\pi_{p^{-1}}\Big(^{p}\big(S_{q^{-1}}(v_{(1,q^{-1})})\cdot w\big)_{(1,prp^{-1})}\Big)\\
 &\stackrel{(\ref{qco})}=&^{p}\Big(\big(S_{q^{-1}}(v_{(1,q^{-1})})\cdot w\big)_{(0)}\Big)\otimes v_{(0)(0)}\\
 &&\otimes v_{(0)(1,r)}\pi_{p^{-1}}\Big(\pi_{p}\big(S_{q^{-1}}(v_{(1,q^{-1})})\cdot w\big)_{(1,p^{-1}prp^{-1}p)}\Big)\\
 &=&^{p}\Big(\big(S_{q^{-1}}(v_{(1,q^{-1})})\cdot w\big)_{(0)}\Big)\otimes v_{(0)(0)}\otimes v_{(0)(1,r)}\big(S_{q^{-1}}(v_{(1,q^{-1})})\cdot w\big)_{(1,r)}\\
 &\stackrel{(\ref{JR2})}=&^{p}\big(S_{q^{-1}}(v_{(1,q^{-1})})_{(2,q)}\cdot w_{(0)}\big)\otimes v_{(0)(0)} \\ 
 &&\otimes v_{(0)(1,r)}\big(S_{q^{-1}}(v_{(1,q^{-1})})_{(3,r)}w_{(1,r)}\big)S^{-1}\pi_{q^{-1}}\big(S_{q^{-1}}(v_{(1,q^{-1})})_{(1,qr^{-1}q^{-1})}\big)\\
 &=&^{p}\big(S_{q^{-1}}(v_{(3,q^{-1})})\cdot w_{(0)}\big)\otimes v_{(0)}\otimes v_{(1,r)}\big(S_{r}(v_{(2,r)})w_{(1,r)}\big)\pi_{q^{-1}}(v_{(4,qr^{-1}q^{-1})})\\
 &=&^{p}\big(S_{q^{-1}}(v_{(1,q^{-1})})\cdot w_{(0)}\big)\otimes v_{(0)}\otimes w_{(1,r)})\pi_{q^{-1}}(v_{(2,qr^{-1}q^{-1})})\\
 &=&C_{V,W}(v_{(0)}\otimes w_{(0)})\otimes w_{(1,r)}\pi_{q^{-1}}(v_{(1,qrq^{-1})})\\
 &=&(C_{V,W}\otimes id)\big(v_{(0)}\otimes w_{(0)}\otimes w_{(1,r)}\pi_{q^{-1}}(v_{(1,qrq^{-1})})\big)\\
 &=&(C_{V,W}\otimes id)\rho_{r}^{V\otimes W}(v\otimes w).
 \end{eqnarray*}

 Thirdly, we can find $C_{V,W}$ satisfies the conditions (4.14) and (4.15). But here we only check the first condition, the other is similar.
 \begin{eqnarray*}
 &&(C_{V,^{W}X}\otimes id_{W})(id_{V}\otimes C_{W,X})(v\otimes w\otimes x)\\
 &\stackrel{(\ref{BZ})}=&(C_{V,^{W}X}\otimes id_{W})(v\otimes {}^{q}(S_{s^{-1}}(w_{(1,s^{-1})})\cdot x)\otimes w_{(0)})\\
 &=&C_{V,^{W}X}\Big(v\otimes {}^{q}\big(S_{s^{-1}}(w_{(1,s^{-1})})\cdot x\big)\Big)\otimes w_{(0)}\\
 &\stackrel{(\ref{BZ})}=&{}^{p}\Big(S_{qs^{-1}q^{-1}}(v_{(1,qs^{-1}q^{-1})})\cdot {}^{q}\big(S_{s^{-1}}(w_{1,s^{-1}})\cdot x\big)\Big)\otimes v_{(0)}\otimes w_{(0)}\\
 &\stackrel{(\ref{HHQ1})}=&{}^{pq}\Big(\pi_{q^{-1}}\big(S_{qs^{-1}q^{-1}}(v_{(1,qs^{-1}q^{-1})})\big)\cdot \big(S_{s^{-1}}(w_{1,s^{-1}})\cdot x\big)\Big)\otimes v_{(0)}\otimes w_{(0)}\\
 &=&{}^{pq}\Big(\pi_{q^{-1}}\big(S_{qs^{-1}q^{-1}}(v_{(1,qs^{-1}q^{-1})})\big)S_{s^{-1}}(w_{1,s^{-1}})\cdot x\Big)\otimes v_{(0)}\otimes w_{(0)}\\
 &=&{}^{pq}\big(S_{s^{-1}}\pi_{q^{-1}}(v_{(1,qs^{-1}q^{-1})})S_{s^{-1}}(w_{1,s^{-1}})\cdot x\big)\otimes v_{(0)}\otimes w_{(0)}\\
 &=&{}^{pq}\Big(S_{s^{-1}}\big(w_{(1,s^{-1})}\pi_{q^{-1}}(v_{(1,qs^{-1}q^{-1})})\big)\cdot x\Big)\otimes v_{(0)}\otimes w_{(0)}\\
 &\stackrel{(\ref{HHQ2})}=&{}^{pq}\big(S_{s^{-1}}(v\otimes w)_{(1,s^{-1})}\cdot x\big)\otimes (v\otimes w)_{(0)}\\
 &=&C_{V\otimes W,X}(v\otimes w, x).
 \end{eqnarray*}

 Finally, we check the condition $C_{^{s}V, ^{s}W}={}^{s}(\cdot)C_{V,W}$. Indeed,
 \begin{eqnarray*}
 C_{{}^{s}V, {}^{s}W}(^{s}v\otimes {}^{s}w)
 &=&{}^{sps^{-1}}\Big(S_{sqs^{-1}}\big(({}^{s}v)_{(1,sqs^{-1}})\big)\cdot {}^{s}w\Big)\otimes ({}^{s}v)_{(0)}\\
 &\stackrel{(\ref{HHQ2})}=&{}^{sps^{-1}}\Big(S_{sqs^{-1}}\big(\pi_{s}(v_{(1,s^{-1}sq^{-1}s^{-1}s)})\big)\cdot {}^{s}w\Big)\otimes {}^{s}(v_{(0)})\\
 &=&{}^{sps^{-1}}\big(\pi_{s}S_{q^{-1}}(v_{(1,q^{-1})})\cdot {}^{s}w\big)\otimes {}^{s}(v_{(0)})\\
 &\stackrel{(\ref{HHQ1})}=&{}^{sps^{-1}}\Big(^{s}\big(\pi_{s^{-1}}\big(\pi_{s}S_{q^{-1}}(v_{(1,q^{-1})})\big)\cdot w\big)\Big)\otimes {}^{s}(v_{(0)})\\
 &=&{}^{sps^{-1}}\Big(^{s}\big(S_{q^{-1}}(v_{(1,q^{-1})})\cdot w\big)\Big)\otimes {}^{s}(v_{(0)})\\
 &=&{}^{sp}\big(S_{q^{-1}}(v_{(1,q^{-1})})\cdot w\big)\otimes {}^{s}(v_{(0)})\\
 &=&{}^{s}(\cdot)\big(^{p}(S_{q^{-1}}(v_{(1,q^{-1})})\cdot w)\otimes v_{(0)}\big)\\
 &\stackrel{(\ref{BZ})}=&{}^{s}(\cdot)C_{V,W}(v\otimes w).
 \end{eqnarray*}
 This completes the proof.$\hfill \square$
 \\

 Similar to (\cite{AFG15}), we can  give the braided $C_{V,W}$ an inverse by the following way.

 \textbf{Proposition \thesection.9} Let $(V, \rho^{V})\in \mathscr Y\mathscr D(H)_{p}$ and $(W, \rho^{W})\in \mathscr Y\mathscr D(H)_{q}$. Then can give the braided $C_{V,W}$ an inverse $C^{-1}_{V, W}$, which is defined by
 \begin{eqnarray*}
 C^{-1}_{V, W}: {}^{V}W\otimes V&\rightarrow& V\otimes W,\\
 C^{-1}_{V, W}({^{p}w\otimes v})&=&v_{(0)}\otimes v_{(1,q)}\cdot w,
 \end{eqnarray*}
 where $p,q\in G$.

 \textbf{\emph{Proof}} For any $v\in V, w\in W$, we have
 \begin{eqnarray*}
 C^{-1}_{V, W}C_{V,W}(v\otimes w)
 &=&C^{-1}_{V, W}({}^{p}(S_{q^{-1}}(v_{(1,q^{-1})}\cdot w))\otimes v_{(0)})\\
 &=&v_{(0)}\otimes v_{(1,q)}\cdot ((S_{q^{-1}}(v_{(2,q^{-1})}))\cdot w)\\
 &=&v_{(0)}\otimes (v_{(1,q)}S_{q^{-1}}(v_{(2,q^{-1})}))\cdot w\\
 &=&v_{(0)}\otimes \epsilon(v_{e})\cdot w\\
 &=&v\otimes w.
 \end{eqnarray*}

 Conversely, for any ${}^{p}w\in {}^{V}W, v\in V$,
 \begin{eqnarray*}
 C_{V, W}C^{-1}_{V,W}({}^{p}w\otimes v)
 &=&C_{V, W}(v_{(0)}\otimes v_{(1,q)}\cdot w)\\
 &=&{}^{p}(S_{q^{-1}}(v_{(1,q^{-1})})\cdot (v_{(2,q)}\cdot w))\otimes v_{(0)}\\
 &=&{}^{p}(S_{q^{-1}}(v_{(1,q^{-1})})\cdot v_{(2,q)})\cdot w)\otimes v_{(0)}\\
 &=&{}^{p}w\otimes v.
 \end{eqnarray*}
 Since then $C_{V,W}$ is an isomorphism with inverse $C^{-1}_{V, W}$.$\hfill \square$
 \\


 As a consequence of the above results, we obtain another main result of this paper.
 \\

 \textbf{Theorem \thesection.10} For a crossed group-cograded Hopf quasigroup $H$, we define $\mathscr Y\mathscr D(H)$ as the disjoint union of all $\mathscr Y\mathscr D(H)_{p}$ with $p\in G$. Then $\mathscr Y\mathscr D(H)$ is a braided crossed category over group $G$.

 \textbf{\emph{Proof}} As for $\mathscr Y\mathscr D(H)$ is a subcategory of the category $\mathscr Y\mathscr D\mathscr Q(H)$, so it is a crossed category.
 Then we only need prove $\mathscr Y\mathscr D(H)$ is braided.

 The braiding in $\mathscr Y\mathscr D(H)$ can be given by Proposition \thesection.8, and the braiding is invertible, its inverse is the family $C^{-1}_{V,W}$, which defined in Proposition \thesection.9.
 Hence it is obvious that $\mathscr Y\mathscr D(H)$ is a braided crossed category. $\hfill \square$

\section*{Acknowledgements}

 The work was partially funded by China Postdoctoral Science Foundation (Grant No.2019M651764), National Natural Science Foundation of China (Grant No. 11601231).

\end {document}